\documentclass[12pt,letterpaper,reqno]{amsart}
\usepackage[T1]{fontenc}
\usepackage{lmodern,amsmath,amssymb,amsthm,mathtools}
\usepackage[a4paper,margin=25mm]{geometry}
\usepackage[hidelinks]{hyperref}
\usepackage{microtype}
\hypersetup{pdftitle={Unique determination of ellipsoids by their intrinsic volumes}}
\newtheorem{theorem}{Theorem}[section]
\newtheorem{lemma}[theorem]{Lemma}

\newtheorem{remark}[theorem]{Remark}
\newcommand{\R}{\mathbb R}
\newcommand{\dd}{\,du}
\title{Unique determination of ellipsoids by their intrinsic volumes}
\author{Sergii Myroshnychenko, Kateryna Tatarko, and Vladyslav Yaskin}
\date{}

\address{S.~Myroshnychenko, Department of Mathematics and Statistics, University of the Fraser Valley, Abbotsford, BC V2S 7M8, Canada}
	
	\email{serhii.myroshnychenko@ufv.ca}

	\address{K.~Tatarko, Department of Pure Mathematics, University of Waterloo, Waterloo, ON N2L 3G1, Canada}
	
	\email{ktatarko@uwaterloo.ca}

	\address{V.~Yaskin, Department of Mathematical and Statistical Sciences, University of Alberta, Edmonton, AB T6G 2N8, Canada}
	
	\email{yaskin@ualberta.ca}

	\thanks{The authors were supported in part by NSERC}

	
\begin{document}

\begin{abstract}
We prove that ellipsoids in $\mathbb{R}^n$, $n\ge 2$, are uniquely determined by their intrinsic volumes up to
congruence, resolving the conjecture posed by Gusakova and Zaporozhets in 2017. 
\end{abstract}

\maketitle

\section{Introduction}
Let $K$ be a convex body in $\mathbb R^n$ and denote by $B_2^n$ the unit Euclidean ball in $\mathbb R^n$. The Steiner formula asserts that for every $t\ge 0$ we have
\[
 \operatorname{Vol}_n(K+tB_2^n)
 =\sum_{k=0}^n\kappa_{n-k}V_k(K)t^{n-k},
\]
where $$\kappa_k=\frac{\pi^{k/2}}{\Gamma(k/2+1)}$$ is the volume of $B_2^k$. 

The coefficients $V_k(K)$, $k=0,1, \ldots, n$, are called the intrinsic volumes. Note that $V_n(K)$ is equal to the volume of $K$, $V_{n-1}(K)$ is a multiple of its surface area, $V_1(K)$ is a multiple of its mean width, and $V_0(K)=1$ for any convex body.

Gusakova and Zaporozhets asked whether ellipsoids in $\mathbb R^n$ are uniquely determined (up to an isometry) by their intrinsic volumes. In the three-dimensional case, this question was answered affirmatively by Petrov and Tarasov
\cite{PT}. 
We solve the problem in all dimensions. 
\begin{theorem}\label{thm:main}
Two ellipsoids in $\R^n$, $n\ge 2$, with the same intrinsic volumes
$V_1,\ldots,V_n$ are congruent.
\end{theorem}
The analogous result for dual volumes was obtained in \cite{MTY}. 

Throughout the paper, we assume that all ellipsoids are full dimensional. The degenerate case can be reduced to subspaces of lower dimensions.

\section{Preliminaries}

To every ellipsoid $E$ whose semiaxes squared are
$\lambda_1,\ldots,\lambda_n>0$, we associate the polynomial
\begin{equation}\label{eq:Q}
 Q(u)=\prod_{i=1}^n(1+\lambda_i u)=\sum_{j=0}^n s_j(Q) u^j.
\end{equation}
Here $s_j=s_j(Q)$ is the $j$th elementary symmetric polynomial in
$\lambda_i$ defined by

$$
s_j = \sum\limits_{1 \leq i_1 < \dots < i_{j} \leq n} \prod\limits_{p = 1}^j \lambda_{i_p},\qquad\mbox{if } j\ge 1,
$$
and $s_0=1$. The following is a consequence of Newton's inequalities  (see \cite[Theorem~53]{HLP}):
\begin{equation}\label{eq:Newton}
    s_j^2 \geq s_{j-1} s_{j+1}, \qquad 1 \le j < n.
\end{equation}

  The coefficients $s_1,\ldots, s_n$  determine the semiaxes of $E$ up to permutation. Our goal is to show that the intrinsic volumes determine these coefficients.

Throughout the paper, we will use the following rescaled version of intrinsic volumes:
$$
v_k(K)=\frac{V_k(K)}{\kappa_k}, \qquad 0\le k\le n.
$$

 Since an ellipsoid is uniquely determined (up to isometry) by the associated polynomial $Q$, it is enough to work with such polynomials and we will use the following notation: 
 $$v_k(Q)=v_k(E).$$
If the ellipsoid or its polynomial $Q$ are clear from the context, we will just write $v_k$.

It will often be convenient to work with a polynomial 
\[
 Q(u)=1+s_1u+\cdots+s_{n-1}u^{n-1}+s_n u^n,
\]
whose  coefficients are positive, but $Q$ is not necessarily associated with an ellipsoid.

For such a polynomial and  $0\le k\le n$,  we define
\begin{align*}
 T_k(u)&=\sum_{p<k}s_pu^p,\\
 A_k(u)&=\sum_{p<k}(k-p)s_pu^p,\\
 B_k(u)&=\sum_{p>k}(p-k)s_pu^p.
\end{align*}
In particular, $A_0=B_n=0.$
Note that $$A_k(u)= kT_k(u)-u T_k'(u)$$
and $$A_{k+1}=A_k +T_{k+1}.$$
Similarly,
$$B_k(u)=u(Q-T_k)' - k (Q-T_k).$$
Whenever we need to emphasize that $A_k$ and $B_k$ are associated with a given polynomial $Q$, we will write $A_k^Q$ and $B_k^Q$.

\section{Auxiliary lemmas}
First we derive a formula for $v_k$ using an existing formula for intrinsic volumes of ellipsoids proved in~\cite{GSZ}.
\begin{lemma}
Let $E$ be an ellipsoid with the corresponding polynomials $Q$ and $T_k$. Then
\begin{equation}\label{eq:W}
v_k=\frac12\int_0^\infty
 u^{-k/2}\frac{T_kQ'-T_k'Q}{Q^{3/2}}\dd,
 \qquad 1\le k \le n.
\end{equation}
    
\end{lemma}

\begin{proof}
 
By \cite[Theorem 1]{GSZ}, we have
\[
 v_k=\sum_{i=1}^n\lambda_i s_{k-1}^{(i)} 
 \int_0^\infty
 \frac{t^{k-1}}
 {(1+\lambda_i t^2)\prod_{\ell=1}^n\sqrt{1+\lambda_\ell t^2}}\,dt,
\]
where $s_{k-1}^{(i)} $ is the $(k-1)$st elementary symmetric function of $\lambda_1, \ldots, \lambda_n$ with $\lambda_i$ omitted.

Set $u=t^2$. Then
$t^{k-1}\,dt=\tfrac12u^{k/2-1}\,du$ and
$\prod_{\ell=1}^n\sqrt{1+\lambda_\ell t^2}=\sqrt{Q(u)}$.
Bringing the sum inside the integral yields
\begin{equation}\label{eq:w vs H}
 v_k=\frac12\int_0^\infty
 \frac{u^{k/2-1}}{\sqrt{Q(u)}}
 \sum_{i=1}^n\frac{\lambda_i s_{k-1}^{(i)}}
                   {1+\lambda_i u}\,du.
\end{equation}

Next we express the sum using only the coefficients of $Q$ and
the partial sum $T_k$. Observe that
\[
 s_p=s_p^{(i)}+\lambda_i s_{p-1}^{(i)},
\]
where  $s_{-1}^{(i)}=0$ by convention.

Thus 
\begin{align*}
 \frac{T_k(u)}{1+\lambda_i u}= &\frac{\sum_{p=0}^{k-1}s_pu^p}{1+\lambda_i u}
 =\frac{\sum_{p=0}^{k-1}s_p^{(i)}u^p+\lambda_i\sum_{p=0}^{k-1}s_{p-1}^{(i)}u^p}{1+\lambda_i u}\\
& =\frac{s_{k-1}^{(i)}u^{k-1}+\sum_{p=0}^{k-2}s_p^{(i)}u^p+\lambda_i\sum_{p=0}^{k-2}s_{p}^{(i)}u^{p+1}}{1+\lambda_i u}=\sum_{p=0}^{k-2}s_p^{(i)}u^p
  +\frac{s_{k-1}^{(i)}u^{k-1}}{1+\lambda_i u}.
\end{align*}
The formula is valid for all $1\le k \le n$. In the case $k=1$, the sum from 0 to $k-2$ is zero.

Now multiply by $\lambda_i$ and sum over $i$. Since 
\begin{equation}\label{eq:derQ/Q}
    \sum_{i=1}^n\frac{\lambda_i}{1+\lambda_i u}=\frac{Q'}Q,
\end{equation}
and
\[
 \sum_{i=1}^n\lambda_i s_p^{(i)}=(p+1)s_{p+1},
\]
we obtain
\[
 \begin{aligned}
 T_k\frac{Q'}Q
 &=\sum_{p=0}^{k-2}(p+1)s_{p+1}u^p
   +u^{k-1}\sum_{i=1}^n
     \frac{\lambda_i s_{k-1}^{(i)}}{1+\lambda_i u}=T_k'+u^{k-1}\sum_{i=1}^n
     \frac{\lambda_i s_{k-1}^{(i)}}{1+\lambda_i u}.
 \end{aligned}
\]
Rearranging,
\[
 \sum_{i=1}^n\frac{\lambda_i s_{k-1}^{(i)}}{1+\lambda_i u}
 =\frac{T_kQ'-T_k'Q}{u^{k-1}Q}.
\]
Substituting this into integral \eqref{eq:w vs H} gives \eqref{eq:W}.

\end{proof}

\begin{remark}
If  $$  Q(u)=\sum_{j=0}^n s_j u^j$$
is a polynomial with positive coefficients, not necessarily associated with an ellipsoid,  we use the right-hand side of \eqref{eq:W}
as the definition of $v_k(Q)$. The identity
\[
 T_kQ'-T_k'Q
 =\sum_{p<k\le q}(q-p)s_ps_q u^{p+q-1}
\]
shows positivity and convergence. The integrand is
$O(u^{k/2-1})$ at zero and $O(u^{(k-n)/2-2})$ at infinity.
\end{remark}

The following lemma is a version of the Alexandrov-Fenchel inequalities for $v_k$.

\begin{lemma} \label{Lem:AF_W}
Let $K\subset \mathbb R^n$ be a convex body. Then
$$ v_k(K)^2\ge v_{k-1}(K)v_{k+1}(K), \qquad 1\le k<n.$$  
\end{lemma}

\begin{proof}
For any convex body $K \subset \mathbb{R}^n$, the intrinsic volumes can be written as
\[
V_k(K)
=\frac1{\kappa_{n-k}}\binom nk V(K[k],B_2^n[n-k]),
\]
 where $V(K[k],B_2^n[n-k])$ is the corresponding mixed volume (see, for example, \cite[(A.26)]{G}).
Thus
\[
v_k(K)=\frac{\binom nk}{\kappa_k\kappa_{n-k}}\,
V(K[k],B_2^n[n-k]).
\]
Dividing by $v_k(B_2^n)$ and using $V(B_2^n[k],B_2^n[n-k])=\kappa_n$  gives
\[
\frac{v_k(K)}{v_k(B_2^n)}
=\frac{V(K[k],B_2^n[n-k])}{\kappa_n}.
\]
Formula (7.66) in \cite{Schneider} shows that this ratio is log-concave
as a sequence in $k$.

Finally observe that the sequence
\[
v_k(B_2^n)
=\frac{\binom nk\kappa_n}{\kappa_k\kappa_{n-k}}
=\frac{\sqrt{\pi}\,\Gamma((n+1)/2)}
       {\Gamma((k+1)/2)\Gamma((n-k+1)/2)}
\]
is also log-concave. Indeed, this follows from the convexity of \(\log\Gamma\), since both arguments of \(\Gamma\) in the denominator are affine functions of \(k\).

 Thus, the sequence $v_k(K)$ is log-concave as the  product of two positive log-concave sequences.
\end{proof}

 In the next lemma, we will derive an alternative representation for $v_k$ and use it in conjunction with the Cauchy-Schwarz inequality to obtain a bound on the ratio of two consecutive intrinsic volumes of ellipsoids.

\begin{lemma}\label{Lem:WNewton}
For any ellipsoid  in $\mathbb R^n$ we have
\begin{equation}\label{eq:adjacent}
 \frac{v_{k+1}}{v_k}<\sqrt{\frac{s_k}{s_{k-1}}},
 \qquad 1\le k<n.
\end{equation}

\end{lemma}
\begin{proof}

Let $\lambda_1,\ldots,\lambda_n>0$ be the squared semiaxes, and write
\[
 Q(u)=\prod_{i=1}^n (1+\lambda_iu)=\sum_{j=0}^n s_j u^j,\qquad
 H_j(u)=\sum_{i=1}^n \frac{\lambda_i s_{j-1}^{(i)}}{1+\lambda_i u}.
\]
Formula \eqref{eq:w vs H} gives
\begin{equation}\label{eq:H-volume}
 v_j=\frac12\int_0^\infty
       \frac{u^{j/2-1}H_j(u)}{\sqrt{Q(u)}}\,du,
 \qquad 1\le j\le n.
\end{equation}

Fix $1\le k<n$ and consider the measure on $(0,\infty)$ given by
\[
 d\mu(u)=\frac{u^{k/2-1}}{2\sqrt{Q(u)}}\,du.
\]
% with positive density.

Then \eqref{eq:H-volume} gives
\[
 v_k=\int_0^\infty H_k(u)\,d\mu(u),\qquad
 v_{k+1}=\int_0^\infty\sqrt u\,H_{k+1}(u)\,d\mu(u).
\]
We will apply the Cauchy--Schwarz inequality to the second integral and use the following two facts:
\begin{equation}\label{eq:H-comparison}
 H_{k+1}(u)\le\frac{s_k}{s_{k-1}}H_k(u),
\end{equation}
and 
\begin{equation}\label{eq:H-moment}
 \int_0^\infty uH_{k+1}(u)\,d\mu(u)=\int_0^\infty H_k(u)\,d\mu(u)=v_k.
\end{equation}

To prove \eqref{eq:H-comparison}, we use the identities
$s_j=s_j^{(i)}+\lambda_i s_{j-1}^{(i)}$, which give
\[
 s_k s_{k-1}^{(i)}-s_{k-1}s_k^{(i)}
 =\lambda_i\left((s_{k-1}^{(i)})^2
                -s_{k-2}^{(i)}s_k^{(i)}\right)\ge0.
\]
The last inequality follows from \eqref{eq:Newton}
%Newton's inequality (see \cite[Theorem~53]{HLP})
for the tuple $\lambda_1,\ldots,\lambda_n$ with $\lambda_i$ removed. Since $s_{-1}^{(i)}=0$ by convention, the inequality also holds for $k=1$. Multiplying by
$\lambda_i/(1+\lambda_i u)$ and summing over $i$, we obtain
\eqref{eq:H-comparison}.

We now prove \eqref{eq:H-moment}.
A direct computation gives
\[
 \begin{aligned}
 uH_{k+1}-H_k
 &=s_k u\sum_{i=1}^n \frac{\lambda_i}{1+\lambda_i u}
       -u\sum_{i=1}^n \frac{\lambda_i^2 s_{k-1}^{(i)}}{1+\lambda_i u}
       -H_k\\
&=s_k \frac{u Q'}{Q}
       -\sum_{i=1}^n  \left(\lambda_i - \frac{\lambda_i}{1+\lambda_i u} \right) s_{k-1}^{(i)}
       -\sum_{i=1}^n \frac{\lambda_i s_{k-1}^{(i)}}{1+\lambda_i u}\\
 &=s_k\left(\frac{uQ'}Q-k\right),
 \end{aligned}
\]
where we used \eqref{eq:derQ/Q},  $u\lambda_i^2/(1+\lambda_i u)
=\lambda_i-\lambda_i/(1+\lambda_i u)$ and
$\sum_{i=1}^n \lambda_i s_{k-1}^{(i)}=k s_k$.
After multiplication by $d\mu$, the last expression can be written as a total
derivative:
\[
 s_k\left(\frac{uQ'}Q-k\right)
       \frac{u^{k/2-1}}{2\sqrt Q}\,du
 =-s_k\left(\frac{u^{k/2}}{\sqrt Q}\right)'du.
\]
Since $0<k<n$, $u^{k/2}/\sqrt Q$ vanishes at the endpoints of $(0, \infty)$. Thus, integrating the preceding equality yields \eqref{eq:H-moment}.

Combining \eqref{eq:H-comparison} and \eqref{eq:H-moment}, we obtain
\begin{align*}
  v_{k+1}^2
 &=\left(\int_0^\infty\sqrt u\,H_{k+1}(u)\,d\mu(u)\right)^2 \\
 & <\left(\int_0^\infty H_{k+1}(u)\,d\mu(u)\right)
    \left(\int_0^\infty uH_{k+1}(u)\,d\mu(u)\right)
 \le\frac{s_k}{s_{k-1}} v_k^2.
\end{align*}
The Cauchy-Schwarz inequality is strict because
$H_{k+1}\,d\mu$ has positive density on $(0,\infty)$, where
$\sqrt u$ is not constant. This proves
\eqref{eq:adjacent}.

All integrals above are finite. Indeed, $H_j(u)=O(1)$ at zero and
$H_j(u)=O(u^{-1})$ at infinity. The density of $H_{k+1}\,d\mu$
is therefore $O(u^{k/2-1})$ at zero and
$O(u^{(k-n)/2-2})$ at infinity. Both its total integral and its 
first moment converge for $1\le k<n$.
\end{proof}

We will also need the following elementary fact about matrix invertibility. 
\begin{lemma}\label{lem:matrix}
Let $M = (M_{k,j})$ be an $m\times m$ matrix with positive entries. Suppose
\[
 M_{k-1,j}-2M_{k,j}+M_{k+1,j}\ge0,\qquad \text{ for }\  1\leq k\ne j \leq m,
\]
where $M_{0,j} = 0 = M_{m+1,j}$ by convention.
Then $M$ is invertible.
\end{lemma}

\begin{proof}
Consider an $m\times m$ matrix $G = (G_{k,j})$ such that
\[
 G_{k,j}=2M_{k,j}-M_{k-1,j}-M_{k+1,j}.
\]
By assumption, $G_{k,j}\le0$ when $k\ne j$ . Also, note that 
\[
 \sum\limits_{k=1}^m G_{k,j}=M_{1,j}+M_{m,j}>0
\]
where we used that $M$ has positive entries. This implies that
\begin{equation}\label{eq:dominance}
G_{j,j}> - \sum\limits_{k\ne j} G_{k,j} =  \sum_{k\ne j}|G_{k,j}| \geq 0.
\end{equation}
We claim that $G$ is invertible. Indeed, assume that there is $x \in \mathbb{R}^m, \ x \ne 0$ such that $G^Tx = 0$, and choose $j$ such that $|x_j| = \max\limits_{1 \leq k \leq m} |x_k| \ne 0$. 
%Note that $|x_j| \ne 0$ since $x \ne 0.$ 
Then, using  \eqref{eq:dominance}, we have
$$
%G_{j,j}>  \sum_{k\ne j}|G_{k,j}| \geq \sum_{k\ne j}|G_{k,j}| \frac{|x_k|}{|x_j|}
G_{j,j}|x_j|
=\left|\sum_{k\ne j}G_{k,j}x_k\right|
\le \sum_{k\ne j}|G_{k,j}|\,|x_k|
\le |x_j|\sum_{k\ne j}|G_{k,j}|< G_{j,j}|x_j|,
$$
which is impossible. Thus, $G$ is invertible. 

Since the rows of $G$ are linear combinations of the rows of $M$, $Mx = 0$ implies $Gx = 0$. Thus, $M$ is also invertible.
\end{proof}

\section{Main Lemmas}
Consider the space of polynomials
$Q(u)=1+s_1u+\cdots+s_{n-1}u^{n-1}+s_n u^n$ with positive coefficients, where the coefficient $s_n$ is fixed. This assumption comes from the fact that we are working with ellipsoids of fixed volume. Thus, this space of polynomials can be identified with $(0, \infty)^{n-1}$, and we can think of  $v_k(Q)$, $1\leq k \leq n-1,$ as a function of $(s_1,\ldots, s_{n-1})\in (0, \infty)^{n-1}$. The lemma below gives the derivatives of this function.

\begin{lemma}\label{lem:derivatives}
For $1\le k,j<n$,
\begin{equation}\label{eq:J}
 \frac{\partial v_k}{\partial s_j}
 =
\displaystyle  \begin{cases}
\displaystyle \frac14 \int_0^\infty u^{j- \frac{k}2-1}Q^{-\frac32} A_k\, du,&\quad j\ge k,\\[1.5ex]
\displaystyle \frac14 \int_0^\infty u^{j- \frac{k}2 -1}Q^{-\frac32} B_k \, du,& \quad j<k.
 \end{cases}
\end{equation}
When $j=k$, either $A_k$ or $B_k$ gives the same integral.

In particular, all the partial derivatives ${\partial v_k}/{\partial s_j}$ are positive for  $1\leq k,j<n$.
\end{lemma}

\begin{proof}
To obtain \eqref{eq:J}, we treat the cases $j\ge k$ and $j < k$ separately.

First, let $j\ge k$. In this case, $T_k$ does not depend on $s_j$. Differentiating~\eqref{eq:W} gives
\begin{align}\label{eq:L3.1_j>k}
    \frac{\partial v_k}{\partial s_j} &= \frac12  \int_0^\infty  u^{-\frac{k}2} \left( (j T_k u^{j-1} - T_k' u^j) Q^{-\frac32} -  \frac32 (T_k Q' - T_k' Q) Q^{-\frac52} u^j \right) \, du \nonumber\\ 
    & = \frac12  \int_0^\infty  u^{j -\frac{k}2 -1}  Q^{-\frac32} \left( j T_k + \frac12 u T_k' - \frac32 u T_k Q' Q^{-1}\right) \, du.
\end{align}
The differentiation is justified since the integrand in \eqref{eq:L3.1_j>k}  is of order $O(u^{k/2-1})$ at zero and $O(u^{(k-n)/2-2})$ at infinity, where the bounds in $O$ are locally uniform in $s_1, \ldots, s_{n-1}$.

Now we observe that
$$
\frac12 \left( u^{j-\frac{k}2}T_kQ^{-\frac32} \right)' = \frac12 \left(j - \frac{k}2 \right)u^{j - \frac{k}2 - 1} T_k Q^{-\frac32} + \frac12 u^{j-\frac{k}2} T_k'  Q^{-3/2} - \frac34 u^{j-\frac{k}2} T_k Q^{-\frac52} Q'.
$$
The boundary values of $ u^{j-k/2}T_kQ^{-3/2}$ vanish at $0$ and at $\infty$ since $j-\frac{k}2>0$ and  $j + \frac{k}2 - 1 < \frac{3n}2$.
Subtracting this total derivative from the integrand in \eqref{eq:L3.1_j>k} yields
$$
\frac{\partial v_k}{\partial s_j} = \frac14\int_0^\infty u^{j -\frac{k}2 -1} Q^{-\frac32} \left(k T_k - u T_k' \right)  \, du = \frac14\int_0^\infty u^{j -\frac{k}2 -1} Q^{-\frac32} A_k \, du
$$
which gives  \eqref{eq:J} for $j \geq k.$

Now, let $j < k.$  We can rewrite~\eqref{eq:W} as 
$$
v_k=\frac12\int_0^\infty
 u^{-\frac{k}2}\frac{(Q - T_k)'Q-(Q - T_k)Q'}{Q^{\frac32}}\dd
$$
where $Q - T_k$ does not depend on $s_j$. Performing a similar computation as in~\eqref{eq:L3.1_j>k} and adding the total derivative $ \frac12 \left( u^{j-k/2}(Q - T_k)Q^{-3/2} \right)'$ to the integrand  gives
$$
\frac{\partial v_k}{\partial s_j} = \frac14\int_0^\infty u^{j -\frac{k}2 -1} Q^{-\frac32} \left(u (Q - T_k)' - k (Q - T_k) \right)  \, du = \frac14\int_0^\infty u^{j -\frac{k}2 -1} Q^{-\frac32} B_k \, du.
$$
The boundary values also vanish at $0$ and at $\infty$. This confirms~\eqref{eq:J} for $j < k.$

Finally, since $A_k-B_k=kQ-uQ'$, we get
\begin{align*}\label{eq:bridge}
 \int_0^\infty u^{\frac{k}2-1}(A_k-B_k)Q^{-\frac32}\, du &= \int_0^\infty u^{\frac{k}2-1} (kQ-uQ') Q^{-\frac32}\, du \nonumber \\
 &= 2 \int_0^\infty \left(\frac{u^{\frac{k}2}}{\sqrt{Q}}  \right)' \, du =0
\end{align*}
where in the last equality the boundary terms vanish at $0$ and at $\infty$ since $0< k< n$.
Thus, when $j = k$, either expression in~\eqref{eq:J} can be used.

The coefficients of $A_k$ and $B_k$ are positive, so $\frac{\partial v_k}{\partial s_j}$ are also positive for any $1\leq k,j<n$.
\end{proof}

To the integrands in~\eqref{eq:J}, we associate the normalized functions
\[
 a_k(u)=\frac{u^{-k/2}A_k(u)}{v_k},\qquad
 b_k(u)=\frac{u^{-k/2}B_k(u)}{v_k},
 \qquad 0\le k \leq n.
\]
Here $v_0=1$, and $a_0=b_{n}=0$. The normalization by $v_k$ will be useful when comparing two ellipsoids with the same intrinsic volumes.

\begin{lemma}\label{lem:convex}
Let $E$ be an ellipsoid with the corresponding polynomial $Q$. For every $u>0$, the
sequences $a_k(u)$ and $b_k(u)$ are strictly convex, that is, 
\[
 a_{k-1}-2a_k+a_{k+1}>0,\qquad
 b_{k-1}-2b_k+b_{k+1}>0
\]
for every $ 1\leq k<n.$
\end{lemma}

\begin{proof}
First, we treat the sequence $a_k(u)$.
Fix $1\leq k<n$.
%and put $r=v_{k+1}/v_k$. 
By Newton's inequalities \eqref{eq:Newton}, we have
\[
 \frac{s_p}{s_{p-1}}\ge\frac{s_k}{s_{k-1}}>\left( \frac{v_{k+1}}{v_k}\right)^2,
 \qquad 1\leq p\le k,
\]
where we use Lemma~\ref{Lem:WNewton}. Therefore
\begin{equation}\label{eq:positive-polynomial}
 T_{k+1}-\left(\frac{v_{k+1}}{v_k}\right)^2 uT_k
 =s_0+\sum_{p=1}^k\left(s_p-\left(\frac{v_{k+1}}{v_k}\right)^2 s_{p-1}\right)u^p>0
\end{equation}
since every coefficient in this expression is positive.

We now consider
\begin{align*}
 u^{\frac{k+1}2}v_{k+1}(a_{k-1}-2a_k+a_{k+1})
 &=A_{k+1}-2\frac{v_{k+1}}{v_k} \sqrt u\,A_k
   +u\frac{v_{k+1}}{v_{k-1}}A_{k-1}\\
 &\geq A_{k+1}-2 \frac{v_{k+1}}{v_k}\sqrt u\,A_k+ \left(\frac{v_{k+1}}{v_k}\right)^2 uA_{k-1}\\
 &=\left(1-\frac{v_{k+1}}{v_k} \sqrt u\right)^2A_k+T_{k+1}-\left(\frac{v_{k+1}}{v_k}\right)^2 uT_k>0
\end{align*}
where we used Lemma~\ref{Lem:AF_W} and the identities $A_{k+1}=A_k+T_{k+1}$ and $A_{k-1}=A_k-T_k$ .
Since the first term is nonnegative and the second is strictly positive by~\eqref{eq:positive-polynomial},
$$
a_{k-1}-2a_k+a_{k+1} > 0.
$$

We now prove the assertion for $b_k$. Without loss of generality, assume that $E$ is centered at the origin. Consider
\[
 Q^*(u)=\frac{u^n}{s_n} Q\left( \frac1u\right)=\prod_{i=1}^n(1+\lambda_i^{-1}u) = \sum\limits_{p = 0}^n \frac{s_{n-p}}{s_n} u^p.
\]
We note that $Q^*$ is the polynomial of the polar ellipsoid $E^*$ of $E$. By \cite[Prop. 4.8]{KZ}, 
$$
v_j(E^*)=\frac{v_{n-j}(E)}{\sqrt{s_n}}.
$$ 
Since the coefficients of $Q^*$ are $s_p^* = \frac{s_{n-p}}{s_n}$, we get
\begin{align*}
    A_{n-k}^{Q^*}\left(\frac1u\right) &= \sum\limits_{p < n - k} (n-k - p) s_p^* u^{-p} =  \frac1{s_n} \sum\limits_{p < n - k} (n-k - p) s_{n - p} u^{-p} \\
    &= \frac1{s_n} \sum\limits_{k < p} ( p -k ) s_{p} u^{p-n} = \frac{u^{-n}}{s_n} B_k^Q(u).
\end{align*}
Therefore, 
\[
 a_{n-k}^{Q^*}\left(\frac1u\right)= \frac{u^{\frac{n-k}2}A_{n-k}^{Q^*}\left(\frac1u\right)}{v_{n-k}(E^*)} =  \frac{u^{\frac{-n - k}2}B_{k}^{Q}\left(u\right)}{\sqrt{s_n}\, v_{k}(E)} =  \frac{u^{-\frac{n}2}}{\sqrt{s_n}}b_k^Q(u).
\]
Thus, 
\begin{align*}
    b^Q_{k-1}(u)-2b^Q_k(u)+b^Q_{k+1}(u) = \sqrt{s_n} u^{\frac{n}2} \left(a^{Q^*}_{n-k+1}(1/u) - 2 a^{Q^*}_{n-k}(1/u) + a^{Q^*}_{n-k-1}(1/u) \right) > 0
\end{align*}
where the result follows from the assertion for $a_k$ applied to $Q^*.$
\end{proof}

\section{Proof of Theorem~\ref{thm:main}}
%subsections 4.1 and 4.2

\begin{proof}
Let $E_1$ and $E_2$ be ellipsoids with the same intrinsic volumes, and $P,Q$ be the associated polynomials:
$$P(u)=\sum_{i=0}^n s_i(P)u^i; \qquad Q(u)=\sum_{i=0}^n s_i(Q)u^i.$$
Denote 
$$v_k(P)=v_k(Q)=v_k, \qquad 1\le k\le n,$$
and observe that $$s_n(P)=s_n(Q)=v_n^2.$$
Let
\[
 Q_t=(1-t)P+tQ,\qquad 0\le t\le 1.
\]
The coefficients of $Q_t$ are given by $s_i(Q_t) = (1-t) s_i(P)+t s_i(Q)$,  which are positive numbers. Thus, formula \eqref{eq:W} defines $v_k(Q_t)$ for all $0\le t\le 1$. 

Consider the $(n-1)\times (n-1)$ matrix $M=(M_{k,j})$, where
\begin{equation}\label{eq:M-def}
 M_{k,j}=\frac1{v_k}\int_0^1
 \frac{\partial v_k(Q_t)}{\partial s_j}\,dt,
 \qquad 1\le k,j< n.
\end{equation}
By Lemma~\ref{lem:derivatives},  all entries of $M$ are positive.

We now verify the condition of Lemma \ref{lem:matrix} for the matrix $M$.
By formula \eqref{eq:J}, 
\begin{equation}\label{eq:der W_k}
\frac{\partial v_k(Q_t)}{\partial s_j}=\frac14 \int_0^\infty u^{j- \frac{k}2-1}Q_t^{-\frac32} C(Q_t) du,
\end{equation}
where $C(Q_t)$ is either $A_k^{Q_t}$ or $B_k^{Q_t}$, depending on whether $j\ge k$ or $j<k$. In both cases, $C(Q_t)$ depends linearly on the coefficients of $Q_t$, and
so $C(Q_t)=(1-t)C(P)+tC(Q)$. 
Thus, we can evaluate the following integral:
\begin{align}\label{eq:average}
 \int_0^1\frac{C(Q_t)}{Q_t^{3/2}}\,dt
& = C(P) \int_0^1\frac{1-t}{((1-t)P+tQ)^{3/2}}\,dt+C(Q) \int_0^1\frac{t}{((1-t)P+tQ)^{3/2}}\,dt \notag\\
 &=\frac{2}{(\sqrt P+\sqrt Q)^2}
 \left(\frac{C(P)}{\sqrt P}+\frac{C(Q)}{\sqrt Q}\right),
\end{align}
where the last equality follows by direct integration.
%we used the following elementary integrals

 % $$   \int_0^1\frac{1-t}{((1-t)x+ty)^{3/2}}\,dt
% =\frac{2}{\sqrt x(\sqrt x+\sqrt y)^2}, 
% $$
% and
% $$
% \int_0^1\frac{t}{((1-t)x+ty)^{3/2}}\,dt
% =\frac{2}{\sqrt y(\sqrt x+\sqrt y)^2},
%$$
%for any positive real numbers $x$ and $y$.

For $j\ge k$, combining  \eqref{eq:M-def}, \eqref{eq:der W_k}, \eqref{eq:average}, 
and the definition of $a_k$ for the polynomials $P$ and $Q$:
$$ a_k^P=\frac{u^{-k/2}A_k^P}{v_k}, \qquad  a_k^Q=\frac{u^{-k/2}A_k^Q}{v_k},$$
 we get
\begin{equation}\label{eq:M}
 M_{k,j}=\frac12\int_0^\infty
 \frac{u^{j-1}}{(\sqrt P+\sqrt Q)^2}
 \left(\frac{a_k^P}{\sqrt P}+\frac{a_k^Q}{\sqrt Q}\right)\dd.
\end{equation}
For $j<k$, the same formula holds with $b_k$ in place of $a_k$.
When $j=k$, either formula is valid by Lemma \ref{lem:derivatives}.

Fix a column $j$, and set $M_{0,j}=M_{{n},j}=0$.
If $k<j$, we can express each of $M_{k-1,j}$, $M_{k,j}$, and $M_{k+1,j}$ using formula \eqref{eq:M} with the corresponding functions $a$. %In particular, for  $k+1=j$, Lemma \ref{lem:derivatives} alows us to use eirther $a$ or $b$. 
Thus
\begin{align*}
 M_{k-1,j}-2M_{k,j}+M_{k+1,j}
 =\frac12\int_0^\infty
 \frac{u^{j-1}}{(\sqrt P+\sqrt Q)^2}
 \bigg(&\frac{a_{k-1}^P-2a_k^P+a_{k+1}^P}{\sqrt P}\\
 &+\frac{a_{k-1}^Q-2a_k^Q+a_{k+1}^Q}{\sqrt Q}\bigg)\dd>0,
\end{align*}
by Lemma~\ref{lem:convex}. If $k>j$, we use the formula with $b_{k-1}$, $b_k$, $b_{k+1}$ in the same
way. The conventions $M_{0,j}=M_{{n},j}=0$ agree with $a_0=0$ and $b_{n}=0$.
Therefore, by
Lemma~\ref{lem:matrix} the matrix $M$ is invertible.

Now let $$\Delta s_j=s_j(Q)-s_j(P).$$ Observe that  $\Delta s_n = 0$ due to $s_n(P) = s_n(Q)$.
Since
$$\frac{d}{dt} v_k(Q_t)  = \sum_{j=1}^{n-1} \frac{\partial v_k(Q_t)}{\partial s_j} \Delta s_j, $$
integrating this identity and using \eqref{eq:M-def}, we obtain
\begin{equation*}\label{eq:secant}
 0=\frac{v_k(Q)-v_k(P)}{v_k}
   =\sum_{j=1}^{n-1}M_{k,j}\Delta s_j.
\end{equation*}

Since $M$ is invertible, the last equality implies $\Delta s_j=0$ for all $1\le j\le n-1$. Recall also that $\Delta s_n= 0$. Thus $P=Q$, and their roots determine the same squared semiaxes, with multiplicities.  

\end{proof}

\bigskip

\noindent{\bf AI Disclaimer.} ChatGPT was used in exploratory discussions and in the preparation of this manuscript.

\end{document}